\documentclass[aps,prl,reprint,onecolumn,showpacs,superscriptaddress,amsmath,amssymb,nofootinbib]{revtex4-2}

\usepackage{graphicx}
\usepackage{amssymb}
\usepackage{amsthm}
\usepackage{amsmath}
\usepackage{bm}
\usepackage{xfrac}

\usepackage{xcolor}
\usepackage{soul}
\usepackage{tikz}
\usepackage{mathtools}
\usepackage{hyperref}

\theoremstyle{definition}

\newcommand{\tcm}{\textcolor{black}}

\newcommand{\beq}{\begin{equation}}
\newcommand{\eeq}{\end{equation}}

\newcommand{\dt}{\partial_t}
\newcommand{\dtsq}{\partial_t^2}

\newcommand{\dtdt}{\frac{\mathrm{d}}{\mathrm{d}t}}

\def\rspace{\mathcal{R}}
\def\map{\bm{\phi}}

\def\spatdom{\Omega_t}
\def\refdom{\Omega_0}
\def\spatdomf{\Omega_t^f}
\def\spatdoms{\Omega_t^s}
\def\refdomf{\Omega_0^f}
\def\refdoms{\Omega_0^s}
\def\intspatdomf{\int_{\spatdomf}}
\def\intrefdoms{\int_{\refdoms}}
\def\dom{\ \mathrm{d}\Omega}

\def\fint{\Gamma_t^{f}}
\def\sint{\Gamma_t^{s}}
\def\sfint{\Gamma_t^{sf}}
\def\sfintref{\Gamma_0^{sf}}
\def\intspatsf{\int_{\sfint}}
\def\intsfintref{\int_{\sfintref}}
\def\dos{\ \mathrm{d}\Gamma}

\def\normf{\bm{n}^f}
\def\norms{\bm{n}^s}
\def\normsf{\bm{n}^{sf}}
\def\normsref{\bm{n}_0^s}

\def\lbr{\left(}
\def\rbr{\right)}

\newcommand{\eps}{\epsilon}
\def\cE{\mathcal{E}}
\newcommand{\dtdtcE}{\frac{\mathrm{d}\cE}{\mathrm{d}t}}
\def\glenergy{\Psi^f}
\def\glenergybulk{\Psi^f_\text{bulk}}
\def\glenergyint{\Psi^f_\text{int}}

\def\velf{\bm{v}}
\def\phase{c}
\def\densf{\rho}
\def\solidsf{\gamma_{sf}}

\def\sigmaf{\bm{\sigma}^f}
\def\nablag{\nabla_{\Gamma}}

\def\dfdci{\frac{\partial F}{\partial \phase_i}}

\def\fluidsfternary{\varsigma_i}
\def\summi{\sum_{i = 1}^3}
\def\dsolidsfdci{\frac{\partial \solidsf}{\partial \phase_i}}
\def\sigmasf{\bm{\sigma}^{sf}}

\def\proj{\bm{P}_{\Gamma}}

\def\disps{\bm{u}}
\def\dispsfixed{\left. \bm{u} \right|_{\bm{X}}}
\def\vels{\dt{\disps}}
\def\denss{\rho_0^s}
\def\wfixed{\left. W \right|_{\bm{X}}}
\def\defgradfixed{\left. \defgrad \right|_{\bm{X}}}
\def\nablafixed{\nabla_{\bm{X}}}
\def\sigmas{\bm{\sigma}^s}
\def\pkone{\bm{P}}
\def\defgrad{\bm{F}}

\graphicspath{{../}}

\begin{document}

\title{\Large Energy dissipation in elasto-capillary fluid-structure interaction systems involving three immiscible fluids}

\author{S. R. Bhopalam} \thanks{Corresponding author}
\email{sbhopala@purdue.edu}
\author{H. Gomez}
\affiliation{School of Mechanical Engineering, Purdue University, West Lafayette, IN 47907, USA.}

\date{\today}

\maketitle

\noindent {\large \bf ABSTRACT}

\noindent Here, we consider the elasto-capillary fluid-structure interaction problem studied in \cite{bhopalam_cmame_2022}, i.e., three immiscible fluids in contact with an elastic solid. In this article, we show that the solution of this fluid-structure interaction problem satisfies an energy dissipation law.

\noindent {\large \bf Governing equations}

\noindent {\bf Preliminaries:} We use the open sets $\spatdom \in \rspace^d$ and $\refdom \in \rspace^d$ respectively, to denote the spatial and referential domains occupied by a continuum body, where $d$ is the number of spatial dimensions. We hereby refer to $\spatdom$ and $\refdom$ as Eulerian and Lagrangian domains, respectively. We assume $\refdom$ to be fixed in time and its points to be parameterized by the reference coordinates $\bm{X}$. We define a function $\map$ as a mapping from the Lagrangian to Eulerian domains at time $t$ as $\map(\cdot, t) : \refdom \longmapsto \spatdom = \map(\refdom, t)$ such that $\bm{X} \longmapsto \bm{x} = \map(\bm{X}, t) \ \forall \bm{X} \in \refdom$ where $\bm{x}$ denotes the coordinates of the spatial domain. We define the referential displacement and referential velocity as $\disps(\bm{X}, t) \coloneq \map(\bm{X}, t) - \bm{X}$ and $\velf \coloneq \dt \map = \dt \disps$, respectively, where the operator $\dt$ denotes partial time differentiation. We define the deformation gradient as $\defgrad \coloneq \frac{\partial \map}{\partial \bm{X}}$ and the Jacobian determinant as $J \coloneq \det \defgrad$. In what follows, we use subscripts in the definition of spatial and time derivatives. For example, the subscript $\bm{X}$ in $\left. \dt \disps \right|_{\bm{X}}$ indicates that the time derivative has been computed by holding $\bm{X}$ fixed. When no subscript is specified in the time derivative, we assume that the derivative has been computed by holding $\bm{x}$ fixed. Similarly, in the context of spatial derivatives, the subscript $\bm{X}$ in $\nabla_{\bm{X}} \disps$, for example, indicates that the spatial derivative has been computed with respect to $\bm{X}$. When no subscript is specified in the spatial derivative, we assume that the derivative has been taken with respect to $\bm{x}$. 

\noindent In the Fluid-Structure Interaction (FSI) problem we consider here, we decompose $\spatdom$ into two open sets, $\spatdomf$ and $\spatdoms$, such that $\spatdom \coloneq \spatdomf \cup \spatdoms$, and $\spatdomf \cap \spatdoms = \emptyset$. Here, $\spatdomf$ and $\spatdoms$ refer to the spatial configurations of the fluid and solid, respectively. We also define a similar decomposition of $\refdom$ into $\refdomf$ and $\refdoms$, which denote the referential configurations of the fluid and solid. In what follows, we denote the fluid-solid interface in the spatial domain as $\sfint \coloneq \partial \spatdomf \cap \partial \spatdoms$.

\noindent {\bf Fluids:} The governing equations of the fluids here constitute the continuity and the linear momentum balance equations written in the Eulerian domain. We use a thermodynamically consistent phase-field model --- a ternary Navier-Stokes-Cahn-Hilliard model \cite{boyer_2006, boyer_2010} to describe the dynamics of the three immiscible fluids. In what follows, we make two assumptions: a) we use constant density and viscosity for the fluids, and b) we use a single velocity field to describe the motion of the fluids. We define the Ginzburg-Landau free energy density of the fluids as $\glenergy = \frac{12}{\eps}\glenergybulk + \glenergyint$, where $\glenergybulk = \summi \frac{\fluidsfternary}{2} \phase_i^2 \lbr 1 - \phase_i \rbr^2$ and $\glenergyint = \summi \frac{3}{8} \eps \fluidsfternary \lvert \nabla \phase_i \rvert^2$ represent the bulk and interfacial components of the free energy density. Here, $\eps$ is the diffuse interface length scale, $\phase_i \in \left[0,1\right]$ for $i=1,2,3$ is the phase field denoting the volume fraction of the $i\textsuperscript{th}$ fluid, and $\fluidsfternary$ for $i=1,2,3$ is the spreading coefficient defined as $\varsigma_1 = \gamma_{12} + \gamma_{13} - \gamma_{23}$, $\varsigma_2 = \gamma_{12} + \gamma_{23} - \gamma_{13}$ and $\varsigma_3 = \gamma_{13} + \gamma_{23} - \gamma_{12}$ where $\gamma_{ij}$ denotes the surface tension at the interface between fluids $i$ and $j$.

\noindent {\bf Solid:} The governing equations for the solid here are given by the linear momentum balance equation in the Lagrangian domain. In what follows, we assume the solid to be homogeneous, isotropic and nonlinearly elastic.

\noindent  We state the strong form of the ternary FSI problem as follows: find $p : \spatdomf \times \lbr 0, \text{T} \right] \longmapsto \rspace$, $\velf : \spatdomf \times \lbr 0, \text{T} \right] \longmapsto \rspace^d$, $\phase_i : \spatdomf \times \lbr 0, \text{T} \right] \longmapsto \rspace$ for $i=1,2,3$, $\mu_i : \spatdomf \times \lbr 0, \text{T} \right] \longmapsto \rspace$ for $i=1,2,3$ and $\disps : \refdoms \times \lbr 0, \text{T} \right] \longmapsto \rspace^d$ such that
\begin{subequations}
    \begin{alignat}{3}
        & \text{\tcm{Continuity equation (fluid)}} \quad && \nabla \cdot \velf = 0 \quad \quad \quad && \text{in} \ \spatdomf \times \lbr 0, \text{T} \right] 
        \label{eqn:continuity_eqn}
        \\
        & \text{\tcm{Momentum equation (fluid)}} \quad && \densf \lbr \partial_t \velf + \velf \cdot \nabla \velf \rbr = \nabla \cdot \sigmaf \quad \quad \quad && \text{in} \ \spatdomf \times \lbr 0, \text{T} \right] 
        \label{eqn:momentum_eqn}
        \\ 
        & \text{\tcm{Phase-field equation (fluid)}} \quad && \dt \phase_i + \velf \cdot \nabla \phase_i = \nabla \cdot \lbr \frac{M}{\fluidsfternary} \nabla \mu_i \right) \quad \quad \quad && \text{for} \ i = 1,2,3 \ \text{in} \ \spatdomf \times \lbr 0, \text{T} \right] 
        \label{eqn:phfield_eqn_ternary}
        \\
        & \text{\tcm{Auxiliary equation (fluid)}} \quad && \mu_i = \frac{\delta \glenergy}{\delta \phase_i} \ + \ \beta \quad \quad \quad && \text{for} \ i = 1,2,3 \ \text{in} \ \spatdomf \times \lbr 0, \text{T} \right] 
        \label{eqn:chempot_defn_ternary}
        \\
        & \text{\tcm{Momentum equation (solid)}} \quad && \denss \left. \dtsq \disps \right|_{\bm{X}} = \nabla_{\bm{X}} \cdot \bm{P} \quad \quad \quad && \text{in} \ \refdoms \times \lbr 0, \text{T} \right] 
        \label{eqn:momentum_eqn_solid}
        \\
        & \text{\tcm{Boundary conditions (fluid)}} \quad && \velf = 0  \quad \quad \quad && \text{on} \ \fint \ \times \lbr 0, \text{T} \right]
        \label{eqn:dirichlet_vel_extbdary}
        \\
        & \text{\tcm{Boundary conditions (fluid)}} \quad && \normf \cdot \nabla \mu_i = 0 \quad \quad \quad && \text{for} \ i = 1,2,3 \ \text{on} \ \fint \ \times \lbr 0, \text{T} \right] 
        \label{eqn:freeflux_chempot_fluidbdary}
        \\
        & \text{\tcm{Boundary conditions (fluid)}} \quad && \normf \cdot \nabla \phase_i = 0 \quad \quad \quad && \text{for} \ i = 1,2,3 \ \text{on} \ \fint \ \times \lbr 0, \text{T} \right] 
        \label{eqn:freeflux_phfield_extbdary}
        \\
        & \text{\tcm{Boundary conditions (solid)}} \quad && \norms \cdot \disps = 0  \quad \quad \quad && \text{on} \ \sint \ \times \lbr 0, \text{T} \right] 
        \label{eqn:dirichlet_disp_solidbdary}
        \\
        & \text{\tcm{Boundary conditions (solid)}} \quad && \bm{t}_e \cdot \sigmas \norms = 0 \quad \quad \quad && \text{for} \ e = 1,..,d-1 \ \text{on} \ \sint \ \times \lbr 0, \text{T} \right] \label{eqn:tractionfree_solidbdary} \\
        & \text{\tcm{Fluid-solid interface conditions}} \quad \quad && \velf - \vels \circ \map^{-1} = 0 \quad \quad \quad && \text{on} \ \sfint \ \times \lbr 0, \text{T} \right] 
        \label{eqn:fluidvel_soliddisp_fsibdary}
        \\
        & \text{\tcm{Fluid-solid interface conditions}} \quad \quad && \sigmaf \normsf - \sigmas \normsf = \nablag \cdot \sigmasf \quad \quad \quad && \text{on} \ \sfint \ \times \lbr 0, \text{T} \right] \label{eqn:tractionbalance_fsibdary} \\
        & \text{\tcm{Fluid-solid interface conditions}} \quad && \normsf \cdot \nabla \mu_i = 0 \quad \quad \quad && \text{for} \ i = 1,2,3 \ \text{on} \ \sfint \ \times \lbr 0, \text{T} \right] 
        \\
        & \text{\tcm{Wettability condition}} \quad \quad && \normsf \cdot \nabla \phase_i = h_i \quad \quad \quad && \text{for} \ i = 1,2,3 \ \text{on} \ \sfint \ \times \lbr 0, \text{T} \right] \\
        & \text{\tcm{Initial condition (fluid)}} \quad && \velf = \velf_0 \quad \quad \quad && \text{in} \ \spatdomf 
        \label{eqn:initcond_fluid}
        \\
        & \text{\tcm{Initial condition (fluid)}} \quad && \phase_i = \phase_{i,0} \quad \quad \quad && \text{for} \ i = 1,2,3 \ \text{in} \ \spatdomf \\
        & \text{\tcm{Initial condition (solid)}} \quad && \disps = \disps_0 \quad \quad \quad && \text{in} \ \refdoms
        \label{eqn:initcond_solid}
    \end{alignat}
    \label{eqn:nse_chn_ternary_2}
\end{subequations}
$\!\!\!$where $\velf$ is the fluid velocity, $\densf$ is the fluid density, $\sigmaf = -p \bm{I} + 2 \eta \nabla^s \velf - \frac{3}{4} \varsigma_1 \nabla \phase_1 \otimes \nabla \phase_1 - \frac{3}{4} \varsigma_2 \nabla \phase_2 \otimes \nabla \phase_2 - \frac{3}{4} \varsigma_3 \nabla \phase_3 \otimes \nabla \phase_3$ is the fluid Cauchy stress tensor, $p$ is the fluid pressure, $\eta$ is the dynamic viscosity of the fluids, $\nabla^s$ is the symmetrization of $\nabla$, $M$ is the mobility coefficient associated with the diffusive flux of the fluids, $\mu_i$ for $i = 1,2,3$ is the chemical potential, $\beta$ is the Lagrange multiplier used to impose the constraint $\summi \phase_i = 1$, $\denss$ is the density of the solid, $\disps$ is the solid displacement, $\pkone$ is the first Piola-Kirchoff stress tensor of the solid, $\normf$ is the unit normal vector at the fluid boundary $\fint$, $\norms$ is the unit normal vector at the solid boundary $\sint$, $\normsf$ is the unit normal vector at the fluid-solid interface $\sfint$, $\bm{t}_e$ is an orthonormal basis of $\rspace^{d-1}$ that is orthogonal to $\norms$ and $h_i$ is the wettability condition \cite{bhopalam_cmame_2022}. In Eq.~\eqref{eqn:tractionfree_solidbdary}, $\sigmas$ is the solid Cauchy stress tensor defined by $\sigmas = J^{-1} \pkone \defgrad^T$. In Eq.~\eqref{eqn:tractionbalance_fsibdary}, $\sigmasf = \solidsf \bm{P}_{\Gamma}$ is the stress tensor accounting for the fluid-solid surface tension at $\sfint$ \cite{bhopalam_cmame_2022}, where $\solidsf$ is the surface energy density at $\sfint$ and $\bm{P}_{\Gamma}$ is the surface projection tensor. Additionally, $\nablag$ is the surface gradient \cite{gross_book_2011, buscaglia_cmame_2011} on $\Gamma_t^{sf}$ defined by $\nabla_{\Gamma} = \bm{P}_{\Gamma} \nabla$. The variational derivative $\frac{\delta \glenergy}{\delta \phase_i}$ in Eq.~\eqref{eqn:chempot_defn_ternary} is defined as $\frac{\delta \glenergy}{\delta \phase_i} = \frac{\partial \glenergy}{\partial \phase_i} - \nabla \cdot \frac{\partial \glenergy}{\partial \nabla \phase_i}$. In Eqs.~\eqref{eqn:initcond_fluid}--\eqref{eqn:initcond_solid}, the solution variables with the subscript $0$ denote the initial conditions.\\

\noindent {\large \bf Energy dissipation relation} \\

\noindent The energy functional of the fluid-structure interaction problem can be given by
\beq
\cE = \intspatdomf \frac{1}{2} \densf \lvert \velf \rvert^2 \dom + \intspatdomf \glenergy \dom + \intspatsf \solidsf \dos + \intrefdoms \frac{1}{2} \denss \lvert \vels \rvert^2 \dom + \intrefdoms W \dom,
\label{eqn:energy_functional_fsi}
\eeq
where $W$ denotes the strain energy density of the solid. The terms on the right hand side of Eq.~\eqref{eqn:energy_functional_fsi} are as follows: first term represents the kinetic energy of the fluids, the second term represents the free energy associated with the mixing of the fluids, the third term represents the energetic contribution of the solid-fluid surface tension, the fourth term represents the kinetic energy of the solid and the fifth term represents the strain energy of the solid. In what follows, we show that the solution variables that satisfy Eq.~\eqref{eqn:nse_chn_ternary_2} satisfy an energy dissipation law. To derive this energy dissipation law, we independently evaluate the time derivative of all the terms in Eq.~\eqref{eqn:energy_functional_fsi} and assemble them eventually.

\noindent Using the Reynolds transport theorem \cite{scovazzi_lecnotes_2007}, we show that
\beq
\begin{aligned}
\dtdt \intspatdomf \glenergy \dom &= \intspatdomf \summi \lbr \frac{12}{\eps} \dfdci \dt \phase_i + \frac{3}{4} \eps \fluidsfternary \nabla \phase_i \cdot \nabla \lbr \dt \phase_i \rbr \rbr \dom + \intspatsf \glenergy \velf \cdot \normf \dos, \\
& = \underbrace{\intspatdomf \summi \lbr \frac{12}{\eps} \dfdci - \frac{3}{4} \eps \fluidsfternary \Delta \phase_i \rbr \dt \phase_i \dom}_{T_1^{GL}} + \underbrace{\intspatsf \glenergy \velf \cdot \normf \dos}_{T_2^{GL}} + \underbrace{\intspatsf \summi \frac{3}{4} \eps \fluidsfternary \dt \phase_i \nabla \phase_i \cdot \normf \dos}_{T_3^{GL}}.
\end{aligned}
\label{eqn:glenergy_dissipation_ternary}
\eeq
In deriving Eq.~\eqref{eqn:glenergy_dissipation_ternary}, we have used Eqs.~\eqref{eqn:dirichlet_vel_extbdary} and ~\eqref{eqn:freeflux_phfield_extbdary}. To derive the second step of Eq.~\eqref{eqn:glenergy_dissipation_ternary}, we also use the divergence theorem. For convenience, we split Eq.~\eqref{eqn:glenergy_dissipation_ternary} into three terms $T_1^{GL}$, $T_2^{GL}$ and $T_3^{GL}$, each of which we evaluate independently. $T_1^{GL}$ can be re-written as,
\beq
\begin{aligned}
T_1^{GL} &= \intspatdomf \summi \lbr \mu_i - \beta \rbr \lbr \frac{M}{\fluidsfternary} \Delta \mu_i - \nabla \cdot \lbr \velf \phase_i \rbr \rbr \dom, \\
&= -\intspatdomf \summi \frac{M}{\fluidsfternary} \lvert \nabla \mu_i \rvert^2 \dom - \intspatdomf \summi \mu \velf \cdot \nabla \phase_i \dom - \beta \intspatdomf \summi \lbr \frac{M}{\fluidsfternary} \Delta \mu_i - \nabla \cdot \lbr \velf \phase_i \rbr \rbr \dom, \\
&= -\intspatdomf \summi \frac{M}{\fluidsfternary} \lvert \nabla \mu_i \rvert^2 \dom - \intspatdomf \summi \mu \velf \cdot \nabla \phase_i \dom,
\end{aligned}
\label{eqn:t1gl_dissipation_ternary}
\eeq
where we substitute for $\lbr \frac{12}{\eps} \dfdci - \frac{3}{4} \eps \fluidsfternary \Delta \phase_i \rbr$ and $\dt \phase_i$ from Eqs.~\eqref{eqn:phfield_eqn_ternary} and Eq.~\eqref{eqn:chempot_defn_ternary}, respectively, in the first step of Eq.~\eqref{eqn:t1gl_dissipation_ternary}. While deriving the second step of Eq.~\eqref{eqn:t1gl_dissipation_ternary}, we subsequently use the divergence theorem, Eqs.~\eqref{eqn:continuity_eqn} and ~\eqref{eqn:freeflux_chempot_fluidbdary}. The last step in Eq.~\eqref{eqn:t1gl_dissipation_ternary} follows from the property $\intspatdomf \summi \lbr \frac{M}{\fluidsfternary} \Delta \mu_i - \nabla \cdot \lbr \velf \phase_i \rbr \rbr = \intspatdomf \summi \dt \phase_i = \intspatdomf \dt \summi \phase_i = 0$. Now $T_2^{GL}$ can be re-written as
\beq
\begin{aligned}
T_2^{GL} &= \intspatdomf \glenergy \nabla \cdot \velf \dom + \intspatdomf \velf \cdot \nabla \glenergy \dom, \\
&= \intspatdomf \velf \cdot \nabla \glenergy \dom,
\label{eqn:t2gl_dissipation_ternary}
\end{aligned}
\eeq
where we use the divergence theorem and product rule in the first step of Eq.~\eqref{eqn:t2gl_dissipation_ternary}. To derive the second step of  Eq.~\eqref{eqn:t2gl_dissipation_ternary}, we use Eq.~\eqref{eqn:continuity_eqn}. 

\noindent We use standard tensor-calculus operations to derive the following identity:
\beq
\begin{aligned}
\summi \frac{3}{4} \eps \fluidsfternary \nabla \cdot \lbr \nabla \phase_i \otimes \nabla \phase_i \rbr&= \summi \frac{3}{4} \eps \fluidsfternary \nabla \phase_i \Delta \phase_i \ + \ \summi \frac{3}{8} \eps \fluidsfternary \nabla \lbr \nabla \phase_i \cdot \nabla \phase_i \rbr, \\
&= \summi \nabla \phase_i\lbr \beta - \mu_i + \frac{12}{\eps} \dfdci \rbr \ + \ \summi \frac{3}{4} \eps \fluidsfternary \nabla \phase_i \cdot \nabla \nabla \phase_i, \\
&= \summi \nabla \lbr \frac{12}{\eps} \dfdci + \frac{3}{4} \eps \fluidsfternary \nabla \phase_i \cdot \nabla \phase_i \rbr \ - \ \summi \mu_i \nabla \phase_i \ + \ \summi \beta \nabla \phase_i, \\
&= \nabla \glenergy \ - \ \summi \mu_i \nabla \phase_i,
\label{eqn:identity}
\end{aligned}
\eeq
where we use Eq.~\eqref{eqn:chempot_defn_ternary} and the product rule to derive second step of Eq.~\eqref{eqn:identity}. To derive the fourth step in Eq.~\eqref{eqn:chempot_defn_ternary}, we use the property $\summi \beta \nabla \phase_i = \beta \nabla \summi \phase_i = 0$. 

\noindent We now substitute Eqs.~\eqref{eqn:t1gl_dissipation_ternary} and \eqref{eqn:t2gl_dissipation_ternary} for $T_1^{GL}$ and $T_2^{GL}$, respectively in Eq.~\eqref{eqn:glenergy_dissipation_ternary} to get
\beq
\begin{aligned}
\dtdt \intspatdomf \glenergy \dom &= -\intspatdomf \summi \frac{M}{\fluidsfternary} \lvert \nabla \mu_i \rvert^2 \dom \ + \ \intspatdomf \velf \cdot \lbr \nabla \glenergy - \summi \mu_i \nabla \phase_i \rbr \dom \ + \ \intspatsf \summi \frac{3}{4} \eps \fluidsfternary \dt \phase_i \nabla \phase_i \cdot \normf \dos, \\
&= -\intspatdomf \summi \frac{M}{\fluidsfternary} \lvert \nabla \mu_i \rvert^2 \dom \ + \ \intspatdomf \summi \frac{3}{4} \eps \fluidsfternary \velf \cdot \nabla \cdot \lbr \nabla \phase_i \otimes \nabla \phase_i \rbr \dom \ + \ \intspatsf \summi \frac{3}{4} \eps \fluidsfternary \dt \phase_i \nabla \phase_i \cdot \normf \dos, \\
&= -\intspatdomf \summi \frac{M}{\fluidsfternary} \lvert \nabla \mu_i \rvert^2 \dom \ - \ \intspatdomf \summi \frac{3}{4} \eps \fluidsfternary \nabla \velf : \lbr \nabla \phase_i \otimes \nabla \phase_i \rbr \dom, \\
& + \ \intspatsf \summi \frac{3}{4} \eps \fluidsfternary \lbr \dt \phase_i + \velf \cdot \nabla \phase_i \rbr \nabla \phase_i \cdot \normf \dos,
\end{aligned}
\label{eqn:glenergy_dissipation_finalized_ternary}
\eeq
where we use the identity from Eq.~\eqref{eqn:identity} in the second step. To derive the third step of Eq.~\eqref{eqn:glenergy_dissipation_finalized_ternary}, we use the divergence theorem and Eq.~\eqref{eqn:freeflux_phfield_extbdary}. Using the Reynolds transport theorem \cite{scovazzi_lecnotes_2007}, we now show that
\beq
\begin{aligned}
\dtdt \intspatdomf \frac{1}{2} \densf \lvert \velf \rvert^2 &= \intspatdomf \densf \dt \velf \cdot \velf \dom \ + \ \intspatsf \frac{1}{2} \densf \lvert \velf \rvert^2 \velf \cdot \normf \dos, \\
&= \intspatdomf \velf \cdot \lbr - \densf \velf \cdot \nabla \velf + \nabla \cdot \sigmaf \rbr \dom \ + \ \intspatsf \frac{1}{2} \densf \lvert \velf \rvert^2 \velf \cdot \normf \dos, \\
&= - \intspatdomf \frac{1}{2} \densf \velf \cdot \nabla \lvert \velf \rvert^2 \dom \ - \ \intspatdomf \nabla \velf : \sigmaf \dom \ + \ \intspatsf \velf \cdot \sigmaf \normf \dos \ + \ \intspatsf \frac{1}{2} \densf \lvert \velf \rvert^2 \velf \cdot \normf \dos, \\
&= - \intspatdomf \nabla \velf : \sigmaf \dom \ + \intspatsf \velf \cdot \sigmaf \normf \dos, \\
&= - \intspatdomf \nabla \velf : \eta \nabla^s \velf \dom \ + \ \intspatsf \summi \frac{3}{4} \eps \fluidsfternary \nabla \velf : \lbr \nabla \phase_i \otimes \nabla \phase_i \rbr \dom \ + \ \intspatsf \velf \cdot \sigmaf \normf \dos,
\end{aligned}
\label{eqn:kefluid_dissipation_ternary}
\eeq
where we use Eq.~\eqref{eqn:dirichlet_vel_extbdary} in the first step. To derive the second step in Eq.~\eqref{eqn:kefluid_dissipation_ternary}, we use Eq.~\eqref{eqn:momentum_eqn}. To derive the third step in  Eq.~\eqref{eqn:kefluid_dissipation_ternary}, we use the divergence theorem, Eq.~\eqref{eqn:dirichlet_vel_extbdary} and the identity $\densf \velf \cdot \lbr \velf \cdot \nabla \velf \rbr = \frac{1}{2} \densf \velf \cdot \nabla \lvert \velf \rvert^2$. To derive the fourth step in Eq.~\eqref{eqn:kefluid_dissipation_ternary}, we use the divergence theorem and Eq.~\eqref{eqn:continuity_eqn}. The fifth step follows by substituting the definition of $\sigmaf$ and by using the identity $\intspatdomf \nabla \velf : p \bm{I} \dom = \intspatdomf \nabla \cdot \velf p \dom = 0$. We now show that
\beq
\begin{aligned}
\dtdt \intrefdoms \lbr \frac{1}{2} \denss \lvert \vels \rvert^2 + W \rbr \dom &= \intrefdoms \lbr \denss \dt \dispsfixed \cdot \dtsq \dispsfixed \ + \ \dt \wfixed \rbr \dom, \\
&= \intrefdoms \bigg( \dt \dispsfixed \cdot \lbr \nablafixed \cdot \bm{P} \rbr \ + \ \dt \wfixed \bigg) \dom, \\
&= \intsfintref \dt \dispsfixed \cdot \bm{P} \normsref \dos \ - \ \intrefdoms \bm{P} : \nablafixed \lbr \dt \disps \rbr \dom \ + \ \intrefdoms \dt \wfixed \dom, \\
&= \intsfintref \dt \dispsfixed \cdot \bm{P} \normsref \dos, \\
&= \intspatsf \velf \cdot \sigmas \norms \dos,
\end{aligned}
\label{eqn:kesolid_dissipation_ternary}
\eeq
where $\Gamma_0^{sf}$ is the fluid-solid interface in the referential configuration, $\normsref$ is the unit normal vector at the fluid-solid interface in the referential configuration  pointing in the direction from solid to fluid.
The second step in Eq.~\eqref{eqn:kesolid_dissipation_ternary} follows from Eq.~\eqref{eqn:momentum_eqn_solid}. To derive the third step in Eq.~\eqref{eqn:kesolid_dissipation_ternary}, we use the divergence theorem and Eqs.~\eqref{eqn:dirichlet_disp_solidbdary} and \eqref{eqn:tractionfree_solidbdary}. The fourth step in Eq.~\eqref{eqn:kesolid_dissipation_ternary} follows from the identity $\dt \wfixed = \frac{\partial W}{\partial \bm{F}} : \dt \defgradfixed = \bm{P} : \nablafixed \lbr \dt \disps \rbr$. The last step in Eq.~\eqref{eqn:kesolid_dissipation_ternary} is written using Eq.~\eqref{eqn:fluidvel_soliddisp_fsibdary} and the push-forward relation between the stress tractions from the solid in the spatial and referential configurations.

\noindent We follow \cite{brummelen_etal_gruyter_2018} to show that
\beq
\begin{aligned}
\dtdt \intspatsf \solidsf \dos &= \intspatsf \bigg(\dt \solidsf + \lbr \velf \cdot \normf \rbr \lbr \normf \cdot \nabla \solidsf \rbr \bigg) \dos \ + \ \intspatsf \lbr \nablag \solidsf \cdot \velf + \solidsf \nablag \cdot \velf \rbr \dos, \\
&= \intspatsf \bigg(\dt \solidsf + \velf \cdot \nabla \solidsf \bigg) \dos + \intspatsf \solidsf \nablag \cdot \velf \dos, \\
&= \intspatsf \summi \dsolidsfdci \lbr \dt \phase_i + \velf \cdot \nabla \phase_i \rbr \dos \ + \ \intspatsf \solidsf \proj \colon \nabla \velf\dom.
\end{aligned}
\label{eqn:solidsf_dissipation_ternary}
\eeq
We derive the second step in Eq.~\eqref{eqn:solidsf_dissipation_ternary} by re-arranging the terms and using the definition of surface gradient. To derive the third step in Eq.~\eqref{eqn:solidsf_dissipation_ternary}, we use the property $\intspatsf \solidsf \nablag \cdot \velf \dos = \intspatsf \solidsf \bm{P}_{\Gamma} \nabla \cdot \velf \dos = \intspatsf \solidsf \bm{P}_{\Gamma} : \nabla \velf \dos$. We now assemble the terms from Eqs.~\eqref{eqn:glenergy_dissipation_finalized_ternary} -- \eqref{eqn:solidsf_dissipation_ternary} and subsequently simplify them to derive the following energy-dissipation relation,
\beq
\begin{aligned}
\dtdtcE &= -\intspatdomf \summi \frac{M}{\fluidsfternary} \lvert \nabla \mu_i \rvert^2 \dom \ - \ \intspatdomf \nabla \velf : \eta \nabla^s \velf \dom \ + \ \intspatsf \summi \frac{3}{4} \eps \fluidsfternary \lbr \dt \phase_i + \velf \cdot \nabla \phase_i \rbr \nabla \phase_i \cdot \normf \dos, \\
& + \ \intspatsf \nablag \cdot \sigmasf \cdot \velf \dos \ + \ \intspatsf \summi \dsolidsfdci \lbr \dt \phase_i + \velf \cdot \nabla \phase_i \rbr \dos \ + \ \intspatsf \solidsf \proj \colon \nabla \velf\dom, \\
& = -\intspatdomf \summi \frac{M}{\fluidsfternary} \lvert \nabla \mu_i \rvert^2 \dom \ - \ \intspatdomf \nabla \velf : \eta \nabla^s \velf \dom \ + \ \intspatsf \summi \lbr \frac{3}{4} \eps \fluidsfternary \nabla \phase_i \cdot \normf + \dsolidsfdci \rbr \lbr \dt \phase_i + \velf \cdot \nabla \phase_i \rbr \dos,
\end{aligned}
\label{eqn:energy_disspn_fsi_ternary}
\eeq
where we use Eq.~\eqref{eqn:tractionbalance_fsibdary} in the first step. To derive the second step of Eq.~\eqref{eqn:energy_disspn_fsi_ternary}, we use the property $\intspatsf \nablag \cdot \sigmasf \cdot \velf \dos = \intspatsf \nablag \cdot \lbr \solidsf \bm{P}_{\Gamma} \rbr \cdot \velf \dos = - \intspatsf \solidsf \proj \colon \nabla \velf\dos + \int_{\partial \sfint} \solidsf \bm{t} \cdot \velf \ \mathrm{d} \lbr \partial \Gamma \rbr \approx - \intspatsf \solidsf \proj \colon \nabla \velf\dos$. If the ternary FSI problem is driven by static wetting, i.e., $\frac{3}{4} \eps \fluidsfternary \nabla \phase_i \cdot \normf + \dsolidsfdci = 0$ for $i = 1, 2, 3$, the energy dissipation relation simplifies to 
\begin{equation*}
    \boxed{\dtdtcE = -\intspatdomf \summi \frac{M}{\fluidsfternary} \lvert \nabla \mu_i \rvert^2 \dom \ - \ \intspatdomf \nabla \velf : \eta \nabla^s \velf \dom}
\end{equation*}

\noindent {\bf Acknowledgements:} This research was supported by the National Science Foundation (Award no. CBET 2012242), with H.G. listed as the principal investigator. The opinions, findings, and conclusions, or recommendations expressed are those of the authors and do not necessarily reflect the views of the National Science Foundation.

\bibliographystyle{apsrev4-2}
\bibliography{main}

\end{document}